\documentclass[%
reprint,
amsmath,amssymb,
aps,
prb,
]{revtex4-1}

\usepackage{graphicx}
\usepackage{dcolumn}
\usepackage{bm}
\usepackage{adjustbox}
\usepackage{hyperref}
\usepackage{calc}
\usepackage{tikz}

\usepackage{amsthm}
\theoremstyle{definition}

\usepackage{wrapfig}

\begin{document}

\preprint{APS/123-QED}



\title{Writhe-Based Polymer Link Classification Using Machine Learning} 

\author{Jack Beda$^1$, Djordje Mihajlovic$^{1,2}$, Kasturi Barkataki$^1$, Davide Michieletto$^{1,3}$}
\affiliation{$^1$ School of Physics and Astronomy, The University of Edinburgh, Peter Guthrie Tait
Road, Edinburgh, EH9 3FD, UK.}
\affiliation{$^2$ School of Mathematics, The University of Edinburgh, Peter Guthrie Tait
Road, Edinburgh, EH9 3FD, UK.}
\affiliation{$^3$ MRC Human Genetics Unit, Institute of Genetics and Cancer, University of Edinburgh, Edinburgh EH4 2XU, UK}

\date{\today}

\begin{abstract}
Unique and rapid classification of knots and links is an open mathematical problem that is relevant to a range of (bio)physical systems, including polymer melts, DNA, and proteins.  In this paper, we explore a data-driven approach to the classification problem of link topology. Extending the framework introduced in Ref. ~\cite{Sleiman} (Sleiman et al, 2024 \textit{Soft Matter}, 20(1), pp.71-78), we show that a feedforward neural network trained on the writhe density matrix classifies thermally equilibrated configurations of the first six prime links with \(97\%\) accuracy. We demonstrate that this accuracy remains high across a range of temperatures and lengths of link components, while rapidly deteriorating with the addition of topology-altering Gaussian noise; a result consistent with the writhe density matrix containing features sensitive to topology. Our results show that neural networks based on the writhe density matrix efficiently classify two-component links, establishing machine learning as a promising tool for rapid classification of more complex link topologies, e.g. Borromean rings and multi-component links, as the computational cost of exact numerical calculation of topological invariants becomes prohibitive. 
\end{abstract}

\maketitle
\section{Introduction}

Knots and links have been integral to the fabric of human technology, underpinning a wide range of early innovations, from textiles and garments to tools, weaponry, nets, and ornamental crafts~\cite{kaaronen2025ties, Hardy}. Human fascination with knots has progressed from practical utility to intellectual curiosity, which evolved into a rich mathematical discipline with applications in biology, polymer physics~\cite{Tubiana, Klotz}, fluid dynamics~\cite{Enciso, Kleckner}, and field theories, both classical and quantum \cite{Faddeev, Kaku, Witten}. Beginning with the pioneering knot tabulation by Peter Guthrie Tait~\cite{Cargill}, the field has since progressed to the tabulation of all distinct knots up to 20 crossings~\cite{Thistlethwaite, burton:LIPIcs.SoCG.2020.25}, amounting to billions of distinct configurations, while an analogous classification for links remains comparatively limited with complete link classifications up to 13 crossings~\cite{Faullin}. Even in cases where each component of a link is individually unknotted, their collective arrangement can produce non-trivial and highly structured entanglement, making links both mathematically richer and more relevant for modelling real-world systems of interacting filaments~\cite{Tubiana}.

An overarching goal of knot theory is to determine the equivalence of any two given knots or links, meaning whether one can be continuously deformed into the other without cutting or cutting and glueing strands. In practice, this equivalence problem is addressed through the construction and evaluation of topological invariants that remain unchanged under such deformations. A wide range of invariants has been developed, from classical polynomial invariants to more recent algebraic and geometric formulations~\cite{Rolfsen, ozsvath2004holomorphic, khovanov2000categorification}; however, there still persists the central challenge of identifying invariants that are computationally efficient and sufficiently discriminating to distinguish large numbers of distinct knots and links~\cite{BarNatan, Dlotko}. 

In parallel, there has also been a sustained effort to understand how the geometric embedding of a knot or link encodes its underlying topology; i.e. while topology concerns equivalence under continuous deformations, many physically relevant systems are realised as specific geometric configurations in 3-dimensional space, where quantities such as curvature, writhe, and spatial conformation carry meaningful information~\cite{Tubiana,Sleiman,Mihajlovic}. One of the most striking examples of the relationship between intrinsic knot topology and 3-dimensional conformations is seen in how DNA molecules travel in gels or sediment in dense fluids, where there is a linear correlation between electrophoretic mobility and sedimentation on DNA topology~\cite{MichielettoPNAS2015,Stas1996,Carlen2013}. Similarly, polymers with different topologies can be sorted by leveraging their equilibrium size in size-modulated nanochannels~\cite{Marenda2017}. These observations lead to a natural hypothesis that geometric embedding reflects intrinsic topology; deepening the theoretical understanding of this connection is not only fundamentally interesting but can lead to better bottom-up designs of molecular shapes starting from the underlying polymer topology.

\begin{figure*}
\includegraphics[width = \textwidth]{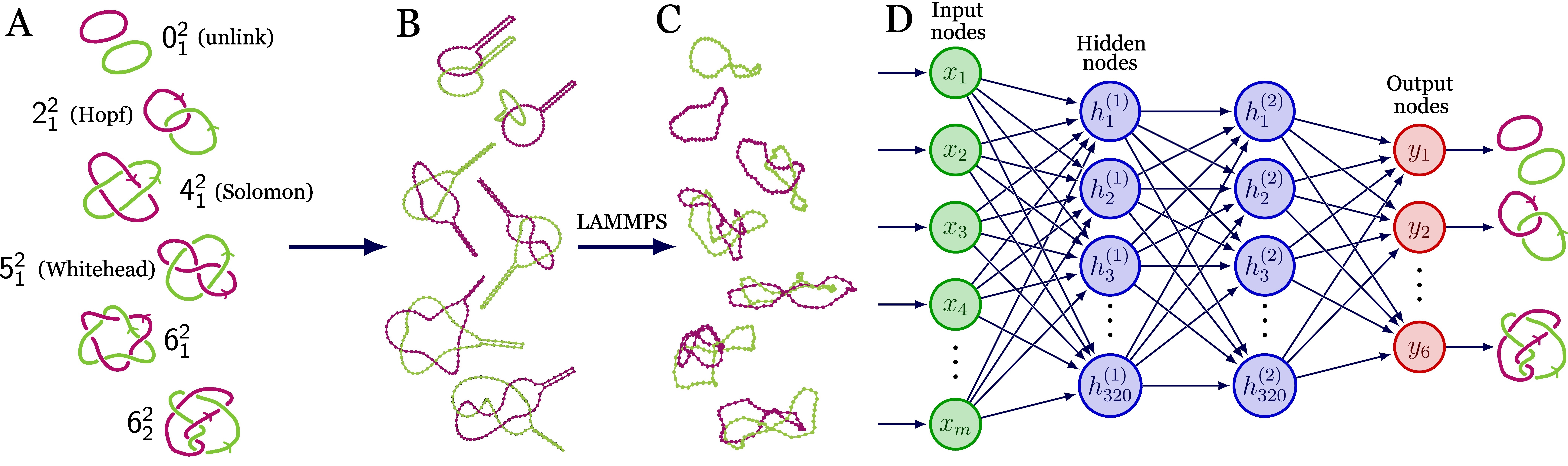}
\caption{\textbf{Overview of the link-generation and classification process}. (A) The first six prime two-component links considered in this work. (B) Corresponding normalized polygonal links used to build the dataset. Each component is extended by a connected sum with a rectangular unknot so that all links have the same number of points per component. (C) Representative thermally agitated configurations generated from the normalized links using LAMMPS. (D) Schematic of the feedforward neural network used for classification. The network is trained on either the flattened \((x,y,z)\) coordinates or the writhe density matrix of a configuration and outputs the predicted link type.}
\label{fig:paper_links_and_lammps_and_NN}
\end{figure*}

The problem of topological classification of knot conformations can be boiled down to pattern recognition. Therefore, many groups have recently considered machine learning to guide the discovery of relationships between existing topological invariants or even the discovery of potentially new invariants~\cite{Davies2021}. Recent work applying machine learning to knot theory has largely focused on diagrammatic inputs, such as braid words or Gauss codes, achieving high classification accuracy but relying on combinatorial rather than geometric information~\cite{Gukov_2021, hughes2020neural}. However, some recent works have also started to explore training deep neural networks (NNs) with geometric data, such as \((x,y,z)\)-coordinates~\cite{Vandans2020,Zhang2026,Braghetto2023} and the writhe density~\cite{Sleiman,Mihajlovic}. Interestingly, training NNs using writhe density data has led to high classification accuracy ($>95\%$) of the first 250 knots with less than 11 crossings~\cite{Sleiman}, a result that was suggested to indicate that NNs may be learning a previously unknown topological invariant~\cite{Sleiman}. While this interpretation remains intriguing, recent findings in Ref.~\cite{Mihajlovic} have identified that the existing NNs used the writhe density to ``shortcut learn'' geometry over topology, i.e. the training datasets could be classified by two ``simple'' elements: (i) the average writhe distributions and (ii) the number of peak and troughs in the writhe density matrix~\cite{Mihajlovic}.  

In view of these results, in this work we extend the framework proposed in Refs.~\cite{Sleiman,Mihajlovic} and use geometric data to build a machine learning classifier of the first six prime links (see Fig.~\ref{fig:paper_links_and_lammps_and_NN}).  On the one hand, our results demonstrate that writhe density matrices provide a rich representation from which NNs can reliably distinguish link types, even under significant variations in discretization and thermal fluctuations. On the other hand, our analysis highlights biases introduced by the constrained exploration of configuration space in molecular dynamics simulations, which can lead to biased training datasets, as highlighted in Ref.~\cite{Mihajlovic}. However, within the context of physical (bio)polymers, our findings suggest that the writhe density matrix encodes geometric features of physically realizable link embeddings that are sufficiently informative for distinguishing link types. In turn, these results highlight our approach as a promising way to build fast link classifiers that can be applied to complex topologies, e.g. polymer melts, Olympic gels, and Borromean rings, where numerical evaluation of exact invariants is often prohibitive.


\section{Methods}

\subsection{Preparation of Links}
\label{sec:links_prep}
In this paper, we consider the first six prime links made of 2 components (see Fig.~\ref{fig:paper_links_and_lammps_and_NN}~A). The links were obtained from KnotPlot~\cite{KnotPlot} as lists of $(x,y,z)$ coordinates and embedded in 3-dimensional space as polygons, i.e. links composed of straight line segments connecting a finite set of vertices. Since feedforward neural networks require fixed input dimensionality, each component must be normalised to have the same number of points. While this can be achieved for knots via interpolation and rescaling, such an approach is not suitable for multi-component links, where components may require different rescalings. Instead, we normalise the links by extending each component through a connected sum with a rectangular unknot of appropriate size. The resulting normalised links are shown in Fig.~\ref{fig:paper_links_and_lammps_and_NN}~B. 

\subsection{LAMMPS Simulation}
\label{sec:lammps_simulation}

For each topological class of link, \(1000\) geometrically distinct configurations were generated by thermally agitating the links using a LAMMPS simulation~\cite{LAMMPS} (see Fig.~\ref{fig:paper_links_and_lammps_and_NN}~C). To simulate the links, each point forming the discretized links are interpreted by LAMMPS as a point mass connected by bonds to its neighbours. Within each component of the link, every point is bonded to its immediate neighbour. LAMMPS uses a Lennard-Jones unit system, based around the Lennard-Jones potential:
\begin{equation}
    U_{\text{LJ}}(r) = 4\epsilon \left[ \left(\frac{\sigma}{r}\right)^{12} - \left(\frac{\sigma}{r}\right)^6 \right],
\end{equation}
where \(\sigma\) and \(\epsilon\) are arbitrary parameters, and \(r\) denotes the radial distance at which the potential is computed. In the Lennard-Jones unit system, the values of the Boltzmann constant (\(k_B\)), \(\sigma\), and \(\epsilon\) are fixed, and all masses are measured relative to a fixed arbitrary mass \(m\). 

The simulations were performed in two stages, namely, an initial equilibration of \(10000\) steps using harmonic bonds, followed by \(200000\) production steps using FENE bonds. The equilibration stage serves to relax the artificially extended link configurations to more natural conformations before introducing the non-linear FENE interactions that preserve topology (see Equations \ref{eq:angular}, \ref{eq:harmonic} and  \ref{eq:fene} below).

Each discretized link is modelled as a chain of particles of mass 1\(m\), evolved with a timestep of \(0.01 \epsilon^{1/2}m^{-1/2}\sigma^{-1}\) under a Langevin thermostat at a temperature of \(1 k_B \epsilon^{-1}\). Non-bonded interactions are described by a Lennard-Jones potential truncated at \(r=2^{1/6}\sigma\), while angles are restricted by the potential 
\begin{equation}\label{eq:angular}
    U_{\text{angular}}(\theta) = K(1+\cos \theta)
\end{equation}
with \(K = 10 \epsilon\). In the equilibration run, bonded particles interact via a harmonic potential
    \begin{equation}\label{eq:harmonic}
        U_{\text{harmonic}}(r) = K (r-r_0)^2
    \end{equation}
with \(K = 100 \epsilon \sigma^{-2}\) and \(r_0 = 1.1 \sigma\). In the production run, bonded particles interact via a FENE potential
    \begin{equation}\label{eq:fene}
        U_{\text{FENE}}(r) = -\frac{1}{2} K R_0^2 \ln \left[ 1 - \left( \frac{r}{R_0} \right)^2 \right]
    + U_{LJ}(r)
    + \epsilon
    \end{equation}
with \(K = 30\epsilon \sigma^{-2}\) and \(R_0 = 1.5 \sigma\) the maximum extent of the bond. Simulations were carried out on \(1000\) copies of each link, producing ensembles of thermally fluctuating configurations that are no longer easily distinguishable by visual inspection (see Fig.~\ref{fig:paper_links_and_lammps_and_NN}~C).

To verify the topology of the data after LAMMPS simulation, we can rely on both computations of the linking number (see Equation \ref{eq:lk_number}) and Jones polynomial~\cite{Rolfsen}. Specifically, we project the embedding of the knot in 3D in an arbitrary direction to access the planar diagram code of the link, and using SageMath~\cite{Sage} to compute the Jones polynomial.

\subsection{The Writhe Density Matrix}
\label{sec_writhe}

The \textit{total writhe} (or simply writhe) of a knot is a quantity that describes the amount of coiling of the knot in 3-dimensional space and assumes real numbers as values.
The total writhe of a knot \(k\) is defined by
\begin{equation}
        \operatorname{Wr}(k) = \frac{1}{4\pi} \int_k \int_k d\vec{r}_1 \times d\vec{r}_2 \cdot \frac{\vec{r}_1 - \vec{r}_2}{|\vec{r}_1- \vec{r}_2|^3}.
        \label{eq:total_writhe}
\end{equation}
where the integrals run over a parameterization of the knot. Note that while the denominator of the integrand vanishes when \(\vec{r}_1 = \vec{r}_2\), the integral remains well-defined~\cite{Calugareanu}. A similar quantity, called the \textit{Gauss linking integral} or simply the \textit{linking number}, can be defined on links with more than one component, i.e. 
\begin{equation}
        \operatorname{Link}_{i,j}(\ell) = \frac{1}{4\pi} \int_{\ell_i} \int_{\ell_j} d\vec{r}_i \times d\vec{r}_j \cdot \frac{\vec{r}_i - \vec{r}_j}{|\vec{r}_i- \vec{r}_j|^3}.
         \label{eq:lk_number}
\end{equation}
where the integrals run over parameterizations of the \(i\)th and \(j\)th components of the link \(\ell\) with \(n\) components.

While the total writhe is neither a knot invariant nor an integer in general, the linking number is always an integer, and is a link invariant~\cite{Rolfsen}. That is to say that two links with different values for the linking number are topologically distinct, and cannot be smoothly deformed into each other. The two concepts of writhe and linking number may be naturally unified in terms of the \emph{writhe density}.   

Given a link \(\ell\) with \(n\) components, where the \(i\)th component is parametrized by \(\vec{r}_i(t)\) for \(t \in I_i \subset \mathbb{R}\), we define the \emph{\(i\)-\(j\) writhe density} of the link as the function
\begin{equation}
    \rho_{i,j} : I_i \times I_j \to \mathbb{R},
\end{equation}
defined by
\begin{equation}
    \rho_{i,j}(t_i, t_j) = \frac{1}{4\pi} \frac{(\dot{\vec{r}}_i(t_i)\times\dot{\vec{r}}_j(t_j)) \cdot ({\vec{r}_i}(t_i)-\vec{r}_j(t_j)) }{|\vec{r}_i(t_i) - \vec{r}_j(t_j)|^3}. \label{eq:writhe_density}
\end{equation}
The writhe density can also be referred to as \emph{segment to segment (StS) writhe}~\cite{Sleiman} and it is clear that the total writhe and linking number may be seen as integrals of the writhe density.

In the case of polygonal links, which are the objects of interest in this work, the writhe density can be discretized to form the \emph{writhe density matrix}. 
Given a polygonal link \(\ell\) with \(n\) components, the writhe density matrix is an \(n\)-by-\(n\) block matrix, whose \(i\)-\(j\)th block consists of the matrix formed by discretizing the \(i\)-\(j\) writhe density function as follows: Suppose the \(i\)th component of the link is formed of \(n_i\) points at positions \(\vec{R}^{(i)}_{p}\) for \(p = 1, \cdots, n_i\). Similarly, component \(j\) is formed of \(n_j\) points at positions \(\vec{R}^{(j)}_q\) for \(q = 1, \cdots, n_j\). Then the \(p\)-\(q\)th entry of the \(i\)-\(j\) block in the writhe density matrix is given by
\begin{equation}
    \frac{1}{4\pi} \frac{(\dot{\vec{R}}^{(i)}_p\times\dot{\vec{R}}_q^{(j)}) \cdot ({\vec{R}_p}^{(i)}-\vec{R}_q^{(j)}) }{|\vec{R}_p^{(i)} - \vec{R}_q^{(j)}|^3}, \label{eq:discrete_writhe_density}
\end{equation}
where \(\dot{\vec{R}}_p^{(i)} = \vec{R}^{(i)}_{p+1 \operatorname{mod} n_i}-\vec{R}^{(i)}_{p}\) and similarly for \(q\). Note that in this study, a highly optimized algorithm from Banchoff~\cite{Banchoff} is used to compute the writhe density matrix as opposed to direct computation of eq.~\ref{eq:discrete_writhe_density}. The writhe density function is symmetric: \({\rho_{i,j}(t_i, t_j) = \rho_{i,j}(t_j, t_i)}\), as is the writhe density matrix. 

Given a link \(\ell\), an approximation of the total writhe of the \(i\)th component is given by the sum of the elements in the \(i\)-\(i\) block of the writhe density matrix (as a discretization of the integral in eq.~\ref{eq:total_writhe}). Similarly, summing the elements in the \(i\)-\(j\)th block for \(i \ne j\) gives an approximation of the linking number between the \(i\)th and \(j\)th component (as a discretization of the integral in eq.~\ref{eq:lk_number}). For the purpose of this paper, the sums of the \(i\)-\(i\) blocks and \(i\)-\(j\) blocks shall henceforth be referred to as the total writhe and linking number, respectively. The sum of the entire matrix shall be called the \emph{total link writhe}.

\begin{figure*}
\centering
\includegraphics[width=\textwidth]{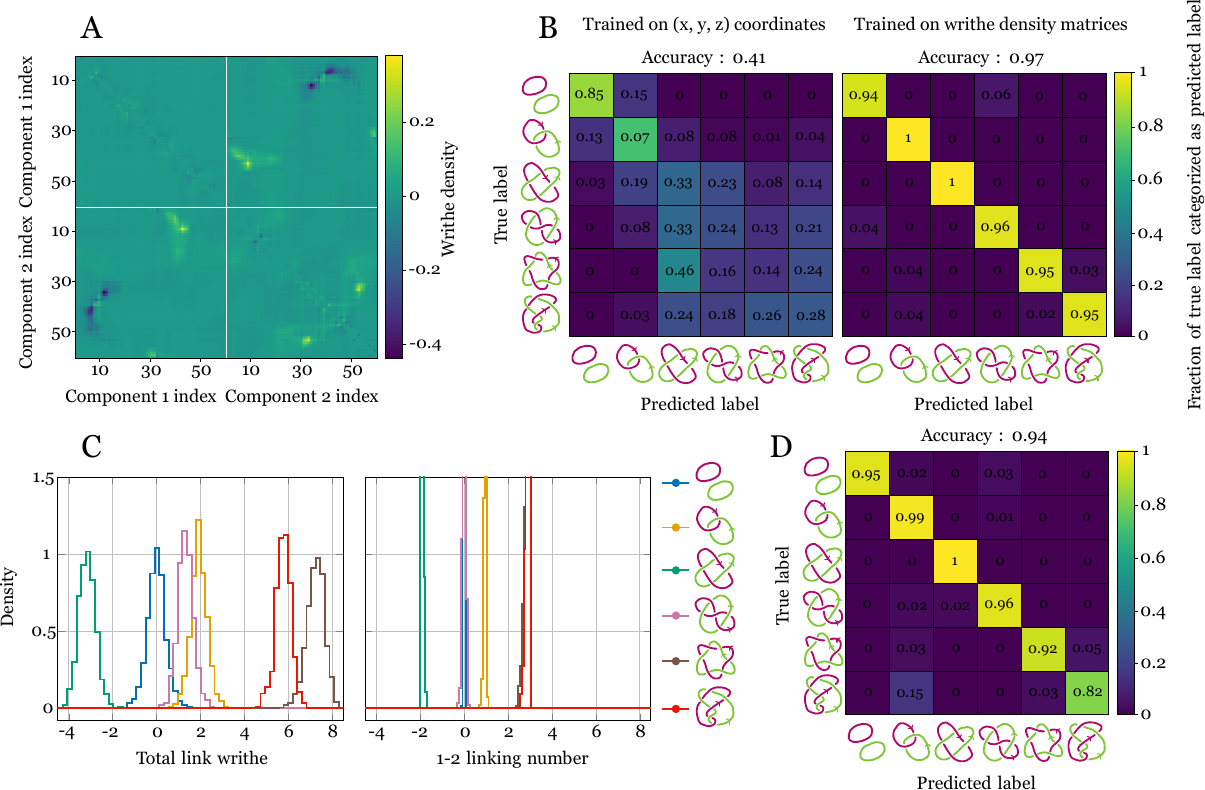}
\caption{\label{fig:paper_writhe_density_and_writhe_distributions_and_confusion_matrices_and_subtraction_of_mean} (A) Representative writhe density matrix for a thermally agitated Whitehead link. The diagonal blocks encode the self-writhe densities of the two components, while the off-diagonal blocks encode their mutual writhe density, whose sum gives the linking number. (B) Confusion matrices for neural networks trained on \((x,y,z)\) coordinates (left) and on full writhe density matrices (right). Training on coordinates gives poor classification performance, whereas training on writhe density matrices yields high accuracy. (C) Distributions of the total link writhe (left) and of the linking number obtained from the off-diagonal block of the writhe density matrix (right) for the six links in the dataset. These global quantities already largely distinguish the link types. (D) Confusion matrix obtained after subtracting the mean of each block of the writhe density matrix so that every block has zero sum, thereby removing the total writhes of the individual components, the total link writhe, and the linking number as features learnable by the neural network. The network still classifies the links with high accuracy, showing that discriminating information remains in the local structure of the writhe density matrix beyond these simple global features.}
\end{figure*}

\subsection{Neural Network Architecture and Training}
\label{sec:training_a_neural_network}

The feedforward neural network (FFNN) used for link classification in this paper consists of three linear layers, each followed by the ReLU activation function  (see Fig.~\ref{fig:paper_links_and_lammps_and_NN}~D). The parameters specifying the feedforward neural network are as follows. The dataset of 6000 links was randomly split into training, validation, and test sets in proportions of 90\%, 7.5\%, and 2.5\%, respectively. The model was trained for 20 epochs using a batch size of 256. Optimization was carried out using the cross-entropy loss function, implemented via \verb|torch.nn.CrossEntropyLoss|, with parameter updates performed using the Adam optimizer (\verb|torch.optim.Adam|)~\cite{Kingma}. Training was conducted within the \verb|pytorch_lightning.Trainer| framework. To mitigate overfitting, early stopping was employed, whereby training was terminated if the validation loss improved by less than 0.05 over 5 consecutive epochs. The code used to build  and train the FFNN is available on the GitLab repository \url{https://git.ecdf.ed.ac.uk/taplab/ml-links-1}. 

\section{Results}

\begin{figure*}
\includegraphics[width = \textwidth]{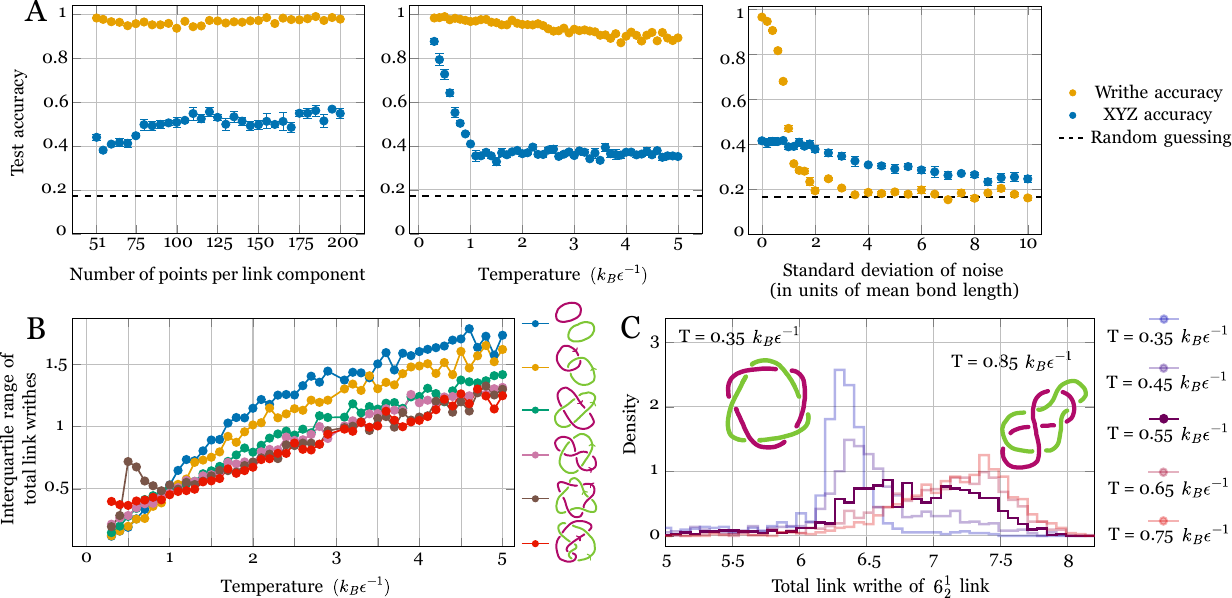}
\vspace{-0.5 cm}
\caption{\textbf{Classification accuracy of writhe-trained FFNNs is robust against temperature and sensitive to Gaussian noice on the coordinates.} (A) Classification accuracy as a function of (left) number of points per link component, (centre) simulation temperature, and (right) addition of Gaussian noise to the $(x,y,z)$ coordinates. (B) Interquartile range of the total link writhe distribution for the six links as a function of simulation temperature. The pronounced jump in the \(6_1^2\) curve near \(T = 0.5 k_B\epsilon^{-1}\) signals a crossover between two preferred configurations. (C) Representative configurations of the \(6_1^2\) link below and above this crossover, illustrating the transition from a rounder and flatter low-temperature geometry to a more elongated  buckled geometry at higher temperatures.}
\vspace{-0.2 cm}
\label{fig:paper_temperatures}
\end{figure*}

\subsection{Link classification is more accurate when FFNNs are trained with the writhe density}

We set off by comparing the performance of FFNNs trained using Cartesian coordinates of the links against networks trained with the writhe density. An example of a writhe density matrix computed on a 2 component link is shown in Fig.~\ref{fig:paper_writhe_density_and_writhe_distributions_and_confusion_matrices_and_subtraction_of_mean}~A). 

In line with Ref. \cite{Sleiman}, the FFNN achieves a low classification accuracy of \(0.41\) when trained on \((x,y,z)\) Cartesian coordinates of the links (see Fig.~\ref{fig:paper_writhe_density_and_writhe_distributions_and_confusion_matrices_and_subtraction_of_mean}~B~left). In contrast, a high accuracy of \(0.97\) is achieved when trained on the writhe density matrices (see Fig.~\ref{fig:paper_writhe_density_and_writhe_distributions_and_confusion_matrices_and_subtraction_of_mean}~B~right). 

At first glance, the success of the FFNN trained on writhe density matrices could be viewed as an indication that the FFNN is able to compute a true topological invariant. However, we also recently discussed the issue of shortcut learning in these problems~\cite{Mihajlovic}; specifically, the distributions of total writhe -- despite not being a topological invariant -- can be seen to be significantly distinct for different link topologies (see Fig.~\ref{fig:paper_writhe_density_and_writhe_distributions_and_confusion_matrices_and_subtraction_of_mean}~C~left). This finding is a consequence of the fact that we realise 3-dimensional physical embeddings of links by performing molecular dynamics simulations of physical polymers, which are subject to volume exclusion, bending rigidity, and inextensible springs between adjacent beads. Removing these constrains leads to much broader writhe distributions that are indistinguishable between different topologies~\cite{Mihajlovic}, however this is not the focus of the current paper.   
The total link writhe is a relatively ``simple'' feature that the NN is likely to learn simply by summing over the elements of the writhe density matrix and constructing relevant bounds for classification. Similarly, the linking number (i.e., the sum of all elements of the \(1\)-\(2\) block of the writhe density matrix), can also be relatively easily learned by the FFNN and separates many of the links (see Fig.~\ref{fig:paper_writhe_density_and_writhe_distributions_and_confusion_matrices_and_subtraction_of_mean}~C~right). Importantly, the linking number is a topological invariant and therefore is expected to be different for different link topologies. Taken together, the distributions of these two quantities effectively separate all the six links considered in this paper and therefore can be used to efficiently train NNs to solve the link classification task. 

However, the performance of the FFNNs can also be examined in the absence of these trivial features by shifting the writhe density matrices of the training data such that all blocks in the writhe density matrix have zero sum. This operation retains the patterns of the writhe density matrices, while removing trivial global features like the total link writhe, the total writhes of the components, and the linking numbers. In terms of Fig.~\ref{fig:paper_writhe_density_and_writhe_distributions_and_confusion_matrices_and_subtraction_of_mean}~C, this operation has the net effect of shifting all distributions to Dirac deltas centred at the origin. After training on the shifted data, the neural network continues to perform very well, with an accuracy of about \(0.94\) (see Fig.~\ref{fig:paper_writhe_density_and_writhe_distributions_and_confusion_matrices_and_subtraction_of_mean}~D). In other words, the writhe density matrix contains distinguishing features beyond the linking numbers and total writhes of the link components, however due to the black box nature of deep NNs, understanding and interpreting these additional features go beyond the scope of this paper.

\subsection{Classification accuracy of writhe-trained FFNNs is robust against longer link lengths and higher simulation temperatures}
\label{sec:variation_with_temperature}

The dataset is created by thermally agitating links made of polymers with fixed number of beads in the LAMMPS molecular dynamics engine using a Langevin thermostat, with temperature $T$ in units of \(k_B \epsilon^{-1}\), and unless otherwise stated, fixed to $T=1 k_B \epsilon^{-1}$. Most of the results shown in this paper are obtained with links made by $N = 60$ points per component. However, we have tested that the accuracy of the network when trained on the writhe density matrices is very robust for longer polymers too. Indeed, the NNs show an accuracy  \(0.96\pm0.01\) across a range of points per component from $N= 51$ to $200$ (see Fig.~\ref{fig:paper_temperatures}~A~left).

 By increasing the value of $T$ in the simulation, we can also explore higher energy configurations of the links and broaden the total link writhe distributions (shown in Fig. \ref{fig:paper_writhe_density_and_writhe_distributions_and_confusion_matrices_and_subtraction_of_mean}~C~left), in turn potentially making the classification task more difficult. The classification accuracy as a function of the simulation temperature is shown in Fig.~\ref{fig:paper_temperatures}~A~(centre). Intriguingly, the accuracy of FFNNs trained on writhe density matrices remain high ($> 90\%$) even at large temperatures. At low temperatures ($T<1 k_B \epsilon^{-1}$), the accuracy of the writhe-trained FFNN is virtually \(100\%\). Across all temperatures, the \((x,y,z)\)-trained FFNNs perform worse than the writhe-trained networks and their accuracy decreases rapidly as $T \simeq 1 k_B\epsilon^{-1}$.  

The robustness of the writhe-trained FFNN accuracy is perhaps surprising. Indeed, as shown in Fig.~\ref{fig:paper_temperatures}~B, as temperature increases, the interquartile range of the total link writhe distributions (i.e. their width) increases, leading to more configurations of topologically distinct links that share the same writhe.  Fig.~\ref{fig:paper_temperatures}~B also shows another interesting feature. The \(6_1^2\) link displays an abrupt change in interquartile width of total link writhe at around \(T=0.5 k_B \epsilon^{-1}\). By inspecting total link writhe distributions at different temperatures we noticed a bimodal, first-order-like transition in the total link writhe (see Fig.~\ref{fig:paper_temperatures}~C). Visual inspection of the \(6_1^2\) link shows that at low temperatures the link  prefers rounder low-writhe conformations,  whereas an elongated configuration is preferred at higher temperatures (see sketches in inset of Fig.~\ref{fig:paper_temperatures}~C). The reason for this buckling may be that at low simulation temperatures, the polymers forming the links are effectively more rigid due to the angular potential modelling persistence length; this increased effective rigidity locks in a round, flat and low-writhe conformation of the \(6_1^2\) link which then buckles at higher temperatures.   

Finally, in an attempt to interpret the accuracy achieved by the FFNNs, we tested the classification of links where Gaussian noise of standard deviation \(\sigma\) was added to each of the \((x,y,z)\) coordinates in the training data. The classification accuracy of the neural network as a function of the noise addition is shown in Fig.~\ref{fig:paper_temperatures}~A~(right). The accuracy of the neural network trained on the writhe density matrices falls sharply, achieving the accuracy of random guesses after the addition of Gaussian noise with standard deviation \(2\). When trained on the \((x,y,z)\) coordinates, however, the accuracy decreases much more slowly, outperforming random guessing until at least \(\sigma = 10\).

It is worth contrasting the response of the neural network to the addition of Gaussian noise (Fig.~\ref{fig:paper_temperatures}~A~right) to the addition of noise from thermal fluctuations (Fig.~\ref{fig:paper_temperatures}~A~centre). In the case of thermal fluctuations, the classification accuracy of a NN trained on writhe density matrices is independent of the temperature of the simulation. Conversely, the accuracy of a NN trained on writhe density matrices plummets as Gaussian noise is added. This is consistent with features learned by the NN present in the writhe density matrices being sensitive to topology, as when the topology is altered by non-topology preserving Gaussian noise, the accuracy decreases abruptly, whereas when the topology is unchanged and thermal noise is added, the accuracy is preserved. On the contrary, the fact that \((x,y,z)\)-trained networks are more robust against Gaussian noise suggest that the learned features are less ``topological'' than those learned by the writhe-trained ones. 

\subsection{NNs learn to distinguish the chirality of the Whitehead Link}

\begin{figure*}
\includegraphics[width = \textwidth]{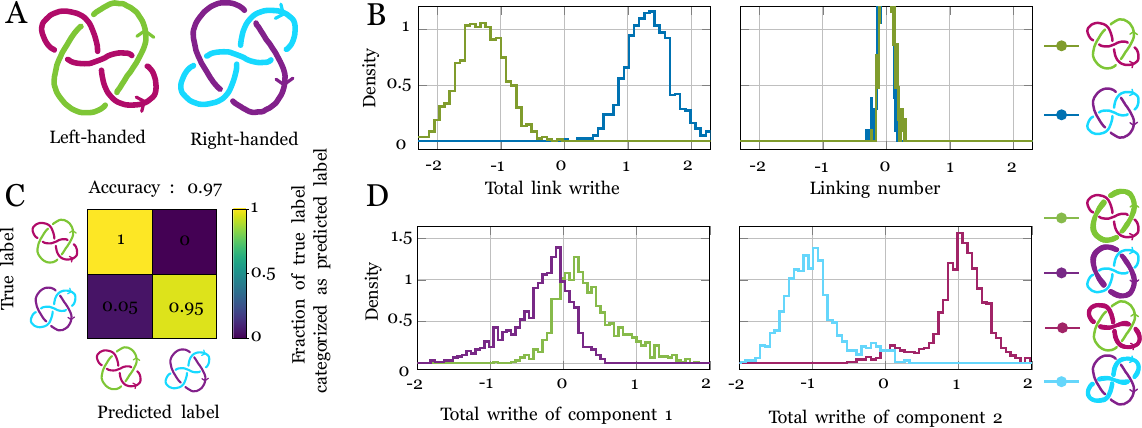}
\caption{ \textbf{NNs learn the chirality of the Whitehead link.} (A) Diagrams of the left-handed Whitehead link and its mirror image. (B) Distributions of the total link writhe (left) and the linking number (right) for thermally agitated configurations of the two chiralities. The total link writhe distinguishes the chiral pair, whereas the linking number distributions are identical. (C) Confusion matrix for a NN trained only on the off-diagonal \(1\)-\(2\) block of the writhe density matrix. Despite the absence of linking number, the network classifies the two chiralities with high accuracy. (D) Distributions of the total writhe of component 1~(left) and component 2~(right) for the two chiralities. The qualitative difference between the two component-wise distributions indicates biased sampling in the simulated ensemble: although the two components of the Whitehead link are topologically equivalent, the initial configuration used in the simulation breaks this symmetry and leads to unequal exploration of configuration space.}
\label{fig:paper_whitehead_case}
\end{figure*}

The Whitehead link is a chiral link, i.e. one which cannot be smoothly deformed to its mirror image (see Fig.~\ref{fig:paper_whitehead_case}~A). The Whitehead link also has linking number \(0\), and so it is possible that the NNs would have difficulty distinguishing the left-handed Whitehead link from its mirror image. To study this, a mirror image of the initial extended Whitehead link from KnotPlot was created by multiplying all coordinates by \(-1\). As before, \(1000\) thermally agitated copies of each chirality were then produced via LAMMPS. When the NN is trained on the full writhe density matrices of the two chiralities of the Whitehead link, it performs exceptionally well (accuracy of \(0.97\), confusion matrix not shown). This is, however, unsurprising, as the distribution of total link writhes of the two chiralities alone can distinguish the chiral pair (see Fig.~\ref{fig:paper_whitehead_case}~B~left).

Since the linking number of both chiralities is \(0\), the distribution of linking numbers of both chiralities are indistinguishable (see Fig.~\ref{fig:paper_whitehead_case}~B~right). Despite this, a NN trained only on the \(1\)-\(2\) block of the writhe density matrices succeeds in classifying the two chiralities with a high accuracy of \(0.97\) (see Fig.~\ref{fig:paper_whitehead_case}~C).

In order to understand this result, we plot the histogram of the total writhes of the first and second components of each of the two chiralities (see Fig.~\ref{fig:paper_whitehead_case}~D). In the figure, the distribution of the total writhe of the first component (Fig.~\ref{fig:paper_whitehead_case}~D~left) differs qualitatively from that of the second component (Fig.~\ref{fig:paper_whitehead_case}~D~right). This indicates that the LAMMPS simulation is not exploring all configurations of the Whitehead link well enough; topologically, both components of the Whitehead link are equivalent (it is possible to deform the link so that the two components switch places~\cite{Skopenkov}). Thus, the initial configuration of the Whitehead link given to LAMMPS (which does not manifestly respect the component interchange symmetry of the Whitehead link) biases the simulation towards exploring some configurations more than others. Again, arguably this is due to the simulated polymers possessing intrinsic rigidity which create an energy barrier to the link exploring some configurations.

Despite evidence of ``biased'' learning, our results are relevant for simulated systems that mimic synthetic and biological polymers, e.g. proteins, DNA, PEG, and others, which do indeed display intrinsic rigidity. For example, if we were able to prepare a Whitehead link made of DNA molecules of length comparable to its persistence length $l_p \simeq 150$ base pairs, we would expect the two components to struggle to switch places as seen in our simulations.

\section{Conclusions}

In this paper, we showed that basic feedforward neural networks can efficiently classify the first six 2-component prime links. We showed that neural networks trained on writhe density matrices achieve classification accuracies of approximately \(97\%\) on thermally agitated configurations of the six links. Consistent with the findings in Ref.~\cite{Sleiman}, the accuracy of writhe-based training of our neural network substantially outperformed networks trained directly on \((x,y,z)\) coordinates.

We examined the robustness of this classification performance under several variations to the training data. The accuracy of the neural network trained on writhe density matrices stayed constant when the number of points discretizing each link component was varied from \(51\) points to \(200\), as well as when the temperature of the LAMMPS simulations was varied from \(0.3 k_B \epsilon^{-1}\) to \(5 k_B \epsilon^{-1}\). 

To investigate whether the network’s success could be attributed solely to simple global features, such as total writhe or linking number, we removed these quantities by shifting the writhe density matrices so that all blocks had zero sum. Even in this setting, the network maintained a high classification accuracy, suggesting that it exploits other structural information in the writhe density matrix. Further evidence of this was provided by the successful classification of the two chiralities of the Whitehead link using only the off-diagonal block of the writhe density matrix, despite the linking number being identical for both chiralities.

When topology-changing Gaussian noise was added to the \((x,y,z)\) coordinates, the accuracy of networks trained on writhe density matrices dropped rapidly to nearly the level of random guessing. Taken together with the temperature-dependence results, this behaviour is consistent with the features learned from the writhe density matrices being sensitive to topology. More precisely, they are preserved under thermal fluctuations that leave the link topology unchanged, but are disrupted by perturbations that alter topology.

Our analysis also highlights important limitations in the generation of training data via molecular dynamics simulations. In particular, for the Whitehead link, we observe biases in the exploration of configuration space, which can lead to systematic distortions in the resulting datasets. This reinforces the conclusion of Ref.~\cite{Mihajlovic} that care must be taken when interpreting strong classification performance as evidence for learned topological invariants. Despite this ``shortcut learning'', our results are relevant for physical (bio)polymer systems where persistence length is an important feature that must be taken into account also when analysing the topology of melts and networks~\cite{He2023}.

In the future, we aim to extend this framework to analyse three component rings which can also form Borromean rings and open curves with the aim of developing computational tools that can rapidly classify complex topologies and can be applied to quantify entanglements in (bio)polymer melts and solutions.  

\section{Data and Code Availability}

\noindent The data and code used to build  and train the FFNN is available on the GitLab repository \url{https://git.ecdf.ed.ac.uk/taplab/ml-links-1}. 

\section*{Acknowledgements}
D.M. acknowledges the Royal Society and the European Research Council (grant agreement No 947918, TAP) for funding. The authors also acknowledge the contribution of the COST Action Eutopia, CA17139. For the purpose of open access, the authors have applied a Creative Commons Attribution (CC BY) licence to any Author Accepted Manuscript version arising from this submission.

\bibliography{bibliography}

\end{document}